\documentclass[12pt,a4paper,reqno]{article}%
\usepackage{amsfonts}
\usepackage{amsmath}
\usepackage{amssymb}
\usepackage{mathrsfs}
\usepackage{graphicx}
\DeclareFontFamily{OT1}{wncyi}{}
\DeclareFontShape{OT1}{wncyi}{m}{it}{
<5> <6> <7> <8> <9> gen * wncyi
<10> <10.95> <12> <14.4> <17.28> <20.74> <24.88> wncyi10
}{}
\DeclareSymbolFont{cyrletters}{OT1}{wncyi}{m}{it}
\DeclareSymbolFontAlphabet{\cyrmath}{cyrletters}
\DeclareMathSymbol{\rE}{\cyrmath}{cyrletters}{003}
\newtheorem{theorem}{Theorem}

\newtheorem{definition}[theorem]{Definition}
\newtheorem{example}[theorem]{Example}

\newtheorem{proposition}[theorem]{Proposition}

\title{\bf Iterated Differential Forms VI:  Differential Equations}
\author{\sc{A.~M.~Vinogradov}\thanks{{\bf e}-{\it mail}: \texttt{vinograd@unisa.it}} and \sc{L.~Vitagliano}\thanks{{\bf e}-{\it mail}: \texttt{luca\_vitagliano@fastwebnet.it}}\\
\small{DMI, Universit\`a degli Studi di Salerno}\\ \small{and INFN, Gruppo collegato di Salerno,}\\
\small{Via Ponte don Melillo, 84084 Fisciano (SA), Italy}}
\begin{document}
\maketitle
\begin{abstract}
We describe the first term of the $\Lambda_{k-1}\mathcal{C}$--spectral
sequence \cite{vv07} of the diffiety $(\mathcal{E},\mathcal{C})$,
$\mathcal{E}$ being the infinite prolongation of an $\ell$--normal system of
partial differential equations, and $\mathcal{C}$ the Cartan distribution on it.
\end{abstract}

\maketitle
\newpage
According to \cite{vv07} secondary $k$--times iterated differential forms and,
in particular, secondary $k$--th degree covariant tensors on a
generic diffiety $(\mathcal{O},\mathcal{C})$ are elements of the
first term of the associated with it
$\Lambda_{k-1}\mathcal{C}$--spectral sequence. In this note we
report results of computation of this term for diffieties that are
infinite prolongations of normal systems of partial differential
equations. To simplify the exposition the equations we deal with
are assumed to be imposed on sections of a vector bundle
$\pi:E\longrightarrow M$, i.e., the diffieties in consideration
are of the form $(\mathcal{E},\mathcal{C})$ with
$\mathcal{E}\subset J^{\infty }(\pi)$ being the infinite
prolongation of a normal system of partial differential equations
$\mathcal{Y}\subset J^s(\pi)$ and $\mathcal{C}$ being the Cartan
distribution on $\mathcal{E}$.

\section{Preliminaries}

Here we collect all necessary preliminaries concerning geometry of
infinitely prolonged PDEs, iterated differential forms (shortly,
IDFs) and secondary IDFs by following
\cite{b99}, \cite{vv06} and \cite{vv07b}, respectively.
Accordingly, we use the notation of these works. As in
\cite{vv07b} we sometimes shorten the notation by using, for
instance, $J^{k}$ instead of  $J^k(\pi)$ for the $k$--th jet
manifold of a fiber bundle $\pi$, etc. As usually $x^1,\dots,x^n$
denote a local chart of $M$.

Let $\xi:V\longrightarrow M$ be a vector bundle, $r=\dim\xi$.
Denote by $v^{1},\ldots,v^{r}$ a system of local linear fiber coordinates
extending $x^1,\dots,x^n$. Any $s$--th order non--linear
differential operator which sends sections of $\pi$ to sections of
$\xi$ may be interpreted as a morphism of bundles
$F:J^{s}\longrightarrow V$. Locally $F$ is represented in the
form
\[
v^{a}=F^{a}(\ldots,x^{\mu},\ldots,u_{\sigma}^{j},\ldots),\quad F^{a}\in
C^{\infty}(J^{k}),\quad a=1,\ldots,r,\,|\sigma|\leq s 
\]
and defines the differential equation $\mathcal{E}_{F}=\{\theta\in
J^{s}\;|\;F(\theta)=0\}\subset J^{s}$. Further on
$\mathcal{E}_{F}$ is assumed to
be regular (see \cite{b99}). Denote by $\mathcal{E}=\mathcal{E}%
_{F}^{\infty}$ the infinite prolongation of $\mathcal{E}_{F}$ and by
$i_{\mathcal{E}%
}:\mathcal{E}\hookrightarrow J^{\infty}$ the corresponding natural
embedding. Put $P=\mathcal{F}\otimes_{C^{\infty}(M)}\Gamma(\xi)$
and recall that there is a unique morphism of diffieties
$\widehat{F}:J^{\infty}\longrightarrow J^{\infty}(\xi)$ such that
$\xi_{\infty,0}\circ\widehat{F}=F\circ\pi _{\infty,k}$.

Let $\Lambda_{k-1}(\xi)$ be the algebra of $(k-1)$--times IDFs on $J^{\infty}(\xi)$ and
$\mathcal{C}_{\circ}\Lambda_{k-1}(\xi)\subset \Lambda_{k-1}(\xi)$
be the $C^{\infty}(V)$--subalgebra generated by elements in the
form $d_{K}^{\,\mathsf{v}}g$, $g\in C^{\infty}(V)\subset C^{\infty
}(J^{\infty}(\xi))$, $K\subset\{1,\ldots,k-1\}$. The pull--back homomorphism $\widehat{F}%
{}^{\ast}:\Lambda_{k-1}(\xi)\longrightarrow\Lambda_{k-1}$ is locally given by
\[
\widehat{F}^{\ast}(d_{L}x^{\mu})=d_{L}x^{\mu},\quad \widehat{F}^{\ast}(d_{K}^{\,\mathsf{v}%
}v_{\sigma}^{a})=(d_{K}^{\,\mathsf{v}}\circ D_{\sigma})(F^{a}),
\]
with $L,K\subset\{1,\ldots,k-1\}$, $\mu=1,\ldots,n=\dim M$,
$a=1,\ldots,r$ and $\sigma$ being a multi--index.
Obviously, $\widehat{F}{}^{\ast}$ sends $\mathcal{C}%
_{\circ}\Lambda_{k-1}(\xi)$ to $\mathcal{C}_{\star}\Lambda_{k-1}$
(see \cite{vv07b}) and this way $\mathcal{C}_{\star}\Lambda_{k-1}$
acquires a structure of a graded
$\mathcal{C}_{\circ}\Lambda_{k-1}(\xi)$--module. By this reason
the notion of a $\mathcal{C}_{\star}\Lambda_{k-1}$--valued
derivation of the algebra $\mathcal{C}_{\circ}\Lambda_{k-1}(\xi)$
is well--defined. By noticing that $C^{\infty
}(M)\subset\mathcal{C}_{\circ}\Lambda_{k-1}(\xi)$ we call a
derivation $Y$ of $\mathcal{C}_{\circ}\Lambda_{k-1}(\xi)$
\emph{vertical} if $Y(C^{\infty}(M))=0$. Denote by
$\Lambda_{k-1}P$ the $\mathcal{C}_{\star}\Lambda_{k-1}$--module of
vertical $\mathcal{C}_{\star}\Lambda_{k-1}$--valued derivations of
$\mathcal{C}_{\circ}\Lambda_{k-1}(\xi)$ and note that $\Lambda
_{0}P\simeq P$.

Let $W_{a}^{K}\in\mathrm{D}(\Lambda_{k-1}(\xi),\Lambda
_{k-1}(\xi))$, $a=1,\ldots,r$, $K\subset\{1,\ldots,k-1\}$, be the
derivations uniquely defined by relations
\begin{align*}
W_{a}^{K}(d_{L}^{\,\mathsf{v}}v_{\sigma}^{b})  &  =\left\{
\begin{array}
[c]{cc}%
1 & \text{if }a=b\text{, }\sigma=\varnothing\text{ and }K=L\\
0 & \text{otherwise}%
\end{array}
\right.  ,\\
W_{a}^{K}(d_L x^{\mu})  &  =0,
\end{align*}
with $b=1,\ldots,r$, $\mu=1,\ldots,n$, $L\subset\{1,\ldots,k-1\}$. Then $W_{a}%
^{K}(\mathcal{C}_{\circ}\Lambda_{k-1}(\xi))\subset\mathcal{C}_{\circ}%
\Lambda_{k-1}(\xi)$, $a=1,\ldots,r$, $K\subset\{1,\ldots,k-1\}$.
Observe that the module $\Lambda_{k-1}P$ is locally generated by
elements $(W_{a}^{K})_{F}=\widehat{F}{}^{\ast}\circ W_{a}^{K}$ and
is horizontal (see \cite{vv07b}). Indeed, the
$C^{\infty}(M)$--module $W_{k-1}(\xi)\subset\Lambda_{k-1}P$
composed of derivations that are locally of the form
$h_{K}^{a}(W_{a}^{K})_{F}$, $h_{K}^{a}\in C^{\infty}(M)$,
$a=1,\ldots,r$, $K\subset\{1,\ldots,k-1\}$ is well defined and
\[
\Lambda_{k-1}P\simeq\mathcal{C}_{\star}\Lambda_{k-1}\otimes_{C^{\infty}%
(M)}W_{k-1}(\xi).
\]

Let $\mathscr{L}\in\mathrm{D}(V)$ be the standard  Liouville
vector field on $V$, i.e.,
$\mathscr{L}=v^{a}\tfrac{\partial}{\partial v^{a}}$. It naturally
extends to a vertical derivation $\mathscr{L}^{\,\mathsf{v}}$ of
$\mathcal{C}_{\circ}\Lambda_{k-1}(\xi)$, which locally looks as
\[
\mathscr{L}^{\,\mathsf{v}}\overset{\mathrm{def}}{=}d_{K}^{\,\mathsf{v}}%
v^{a}\,W_{a}^{K}.
\]
If $\Lambda_{k-1}F\overset{\mathrm{def}}{=}\widehat{F}{}^{\ast}%
\circ\mathscr{L}^{\,\mathsf{v}}\in\Lambda_{k-1}P$, then
$\Lambda_{k-1}F$ has the following local expression
\[
\Lambda_{k-1}F=d_{K}^{\,\mathsf{v}}F^{a}\,(W_{a}^{K})_{F}.
\]

Put $\Lambda_{k-1}\mu_{\mathcal{E}}\overset{\mathrm{def}}{=}\ker
i_{\mathcal{E}}^{\ast}\subset\Lambda_{k-1}$, $i_{\mathcal{E}}^{\ast}%
:\Lambda_{k-1}\longrightarrow\Lambda_{k-1}(\mathcal{E})$, $\mathcal{C}_{\star
}\Lambda_{k-1}\mu_{\mathcal{E}}\overset{\mathrm{def}}{=}$ $\Lambda_{k-1}%
\mu_{\mathcal{E}}\cap\mathcal{C}_{\star}\Lambda_{k-1}$ and $\mathcal{C}%
_{\star}\Lambda_{k-1}(\mathcal{E})\overset{\mathrm{def}}{=}i_{\mathcal{E}%
}^{\ast}(\mathcal{C}_{\star}\Lambda_{k-1})\simeq\mathcal{C}_{\star}%
\Lambda_{k-1}/\mathcal{C}_{\star}\Lambda_{k-1}\mu_{\mathcal{E}}$.
The restriction of a
$\mathcal{C}_{\star}\Lambda_{k-1}$--module $Q$ to $\mathcal{E}$ is defined
to be $Q|_{\mathcal{E}}%
\overset{\mathrm{def}}{=}Q/\mathcal{C}_{\star}\Lambda_{k-1}\mu_{\mathcal{E}%
}\cdot Q$ and $R_{\mathcal{E}}%
^{Q}:Q\ni q\longmapsto\lbrack q]\in Q|_{\mathcal{E}}$ is called
the \emph{restriction homomorphism}.

\section{$\Lambda_{k-1}\mathcal{C}$--modules}

Total derivatives $D_{\mu}, \;\mu=1,\ldots,n $, restricted to
$\mathcal{E}$ extend canonically to derivations of the algebra
$\mathcal{C}_{\star}\Lambda_{k-1}(\mathcal{E})$. These extensions
will be still denoted by $D_{\mu}$'s.

In what follows every
$\mathcal{C}_{\star}\Lambda_{k-1}(\mathcal{E})$--module is
assumed to be graded, locally free and of finite rank. Let
$Q_{1},Q_{2}$ be such ones.

\begin{definition}
A linear differential operator $\Delta:Q_{1}\longrightarrow Q_{2}$
is called $\mathcal{C}$--differential if for any local basis
$\{e_{1},\ldots,e_{t}\}$ of $Q_{1}$, $\Delta$ is locally of the
form $\Delta(p)=(-1)^{|p|\cdot|\alpha
|}\Delta_{\alpha}^{\sigma}D_{\sigma}(p^{\alpha})$ with
$\Delta_{\alpha }^{\sigma}\in Q_{2}$, \;$p=p^{\alpha}e_{\alpha}\in
Q_{1}, \;|\alpha|=
|e_{\alpha}|$ and $p^{\alpha}\in\mathcal{C}_{\star}\Lambda_{k-1}%
(\mathcal{E})$, $\alpha=1,\ldots,t$.
\end{definition}

The left multiplication transforms the totality of all
$\mathcal{C}$--differential operators $\Delta
:Q_{1}\longrightarrow Q_{2}$ into a $\mathcal{C}_{\star}\Lambda
_{k-1}(\mathcal{E})$--module denoted by
$\mathcal{C}\mathrm{Diff}(Q_{1},Q_{2})$. Similarly, the
$\mathcal{C}_{\star}\Lambda _{k-1}(\mathcal{E})$--module of
multiplicity $l, \;Q_{2}$--valued,
multi--$\mathcal{C}$--differential operators on $Q_{1}$, denoted
by $\mathcal{C} \mathrm{Diff}_{(l)}(Q_{1};Q_{2})$, is defined. For
$l=0$ we put
$\mathcal{C}\mathrm{Diff}_{(0)}(Q_{1};Q_{2})=
Q_{2}$. The sub--module of
$\mathcal{C}\mathrm{Diff}_{(l)}(Q_{1},;Q_{2})$,
composed of skew--symmetric operators, is denoted by
$\mathcal{C}\mathrm{Diff}_{(l)}^{\mathrm{alt}}(Q_{1};Q_{2})$ and
\[
\mathrm{alt}_{l}:\mathcal{C}\mathrm{Diff}_{(l)}(Q_{1},Q_{2})\longrightarrow
\mathcal{C}\mathrm{Diff}_{(l)}^{\mathrm{alt}}(Q_{1},Q_{2}).
\]
stands for the alternation operator.

Let $P_{1},P_{2}$ be $\mathcal{C}_{\star}\Lambda_{k-1}%
$--modules and $\square:P_{1}\longrightarrow P_{2}$ a $\mathcal{C}%
$--differential operator.

\begin{proposition}
\label{PropRestr}$\square$ restricts to $\mathcal{E}$, i.e., there exists a
unique $\square|_{\mathcal{E}}\in\mathcal{C}\mathrm{Diff}(P_{1}|_{\mathcal{E}%
},P_{2}|_{\mathcal{E}})$ such that $R_{\mathcal{E}}^{P_{2}}\circ
\square=\square|_{\mathcal{E}}\circ R_{\mathcal{E}}^{P_{1}}$.
\end{proposition}

The module $\mathcal{C}\mathrm{Diff}(\mathcal{C}_{\star}\Lambda
_{k-1}(\mathcal{E}),\mathcal{C}_{\star}\Lambda_{k-1}(\mathcal{E}))$
becomes a unitary ring with
respect to the composition operation. So, the ring
$\mathcal{C}_{\star }\Lambda_{k-1}(\mathcal{E})$ may be viewed as
 a subring of $\mathcal{C}\mathrm{Diff}(\mathcal{C}_{\star
}\Lambda_{k-1}(\mathcal{E}),\mathcal{C}_{\star}\Lambda_{k-1}(\mathcal{E}))$
containing the identity operator. Moreover, if $X\in\mathcal{C}\mathrm{D}(\Lambda_{k-1}%
(\mathcal{E}))$, then $X(\mathcal{C}_{\star}\Lambda_{k-1}(\mathcal{E}%
))\subset\mathcal{C}_{\star}\Lambda_{k-1}(\mathcal{E})$. So, the
restriction
$X|_{\mathcal{C}_{\star}\Lambda_{k-1}(\mathcal{E})}\in\mathcal{C}%
\mathrm{Diff}(\mathcal{C}_{\star}\Lambda_{k-1}(\mathcal{E}),\mathcal{C}%
_{\star}\Lambda_{k-1}(\mathcal{E}))$ is well-defined.

Put
\[
\mathcal{C}\mathrm{D}(\mathcal{C}_{\star}\Lambda_{k-1}(\mathcal{E}%
))\overset{\mathrm{def}}{=}\{X|_{\mathcal{C}_{\star}\Lambda_{k-1}%
(\mathcal{E})}\;|\;X\in\mathcal{C}\mathrm{D}(\Lambda_{k-1}(\mathcal{E}%
))\}\subset\mathcal{C}\mathrm{Diff}(\mathcal{C}_{\star}\Lambda_{k-1}%
(\mathcal{E}),\mathcal{C}_{\star}\Lambda_{k-1}(\mathcal{E})).
\]
As a ring $\mathcal{C}\mathrm{Diff}(\mathcal{C}_{\star}\Lambda_{k-1}(\mathcal{E}%
),\mathcal{C}_{\star}\Lambda_{k-1}(\mathcal{E}))$ is generated by
its submodules
$\mathcal{C}_{\star}\Lambda_{k-1}(\mathcal{E})$ and $\mathcal{C}%
\mathrm{D}(\mathcal{C}_{\star}\Lambda_{k-1}(\mathcal{E}))$.

\begin{definition}
A couple composed of a left graded
$\mathcal{C}_{\star}\Lambda_{k-1}(\mathcal{E})$--module $Q$ and a
homomorphism
\[
\mathcal{C}\mathrm{Diff}(\mathcal{C}_{\star}\Lambda_{k-1}(\mathcal{E}%
),\mathcal{C}_{\star}\Lambda_{k-1}(\mathcal{E}))\ni\Delta\longmapsto\Delta
^{Q}\in\mathcal{C}\mathrm{Diff}(Q,Q)
\]
of unitary rings is called a
$\Lambda_{k-1}\mathcal{C}$\emph{--module} (see, e.g.,
\cite{v01,kv98}).
\end{definition}

\begin{example}
$\mathcal{C}_{\star}\Lambda_{k}^{p}(\mathcal{E})=\mathcal{C}_{\star}%
\Lambda_{k}(\mathcal{E})\cap\Lambda_{k}^{p}(\mathcal{E})\subset\mathcal{C}%
^{p}\Lambda_{k}^{p}(\mathcal{E})$ is a $\Lambda_{k-1}\mathcal{C}$--module for
any $p$. Indeed, for $Y=X|_{\mathcal{C}_{\star}\Lambda_{k-1}(\mathcal{E})}%
\in\mathcal{C}\mathrm{D}(\mathcal{C}_{\star}\Lambda_{k-1}(\mathcal{E}))$,
$X\in\mathcal{C}\mathrm{D}(\Lambda_{k-1}(\mathcal{E}))$,  define
the operator $Y^Q$ by putting
\[
Y^{Q}(\Omega)=\mathcal{L}_{X}^{\{k\}}(\Omega)\in\mathcal{C}_{\star}\Lambda
_{k}^{p}(\mathcal{E}), \quad \Omega
\in\mathcal{C}^{p}\Lambda_{k}^{p}(\mathcal{E}).
\]
Then, in particular, $(\omega Y)^{Q}=\omega Y^{Q}$, $\omega\in\mathcal{C}%
_{\star}\Lambda_{k-1}(\mathcal{E})$, $Y\in\mathcal{C}\mathrm{D}(\mathcal{C}%
_{\star}\Lambda_{k-1}(\mathcal{E}))$.
\end{example}

Note that the tensor product of $\Lambda_{k-1}\mathcal{C}$--modules is a
$\Lambda_{k-1}\mathcal{C}$--module in a natural way.

Let $Q$ be a $\Lambda_{k-1}\mathcal{C}$--module and $P_{1},P_{2}$ be
$\mathcal{C}_{\star}\Lambda_{k-1}(\mathcal{E})$--modules. For $\Delta
\in\mathcal{C}\mathrm{Diff}(P_{1},P_{2})$, $p\in P_{1}$ and $\phi\in
P_{2}^{\ast}=\mathrm{Hom}_{\mathcal{C}_{\star}\Lambda_{k-1}(\mathcal{E}%
)}(P_{2},\mathcal{C}_{\star}\Lambda_{k-1}(\mathcal{E}))$ define $\Delta
(p,\phi)\in\mathcal{C}\mathrm{Diff}(\mathcal{C}_{\star}\Lambda_{k-1}%
(\mathcal{E}),\mathcal{C}_{\star}\Lambda_{k-1}(\mathcal{E}))$ by putting
\[
\Delta(p,\phi)(\omega)\overset{\mathrm{def}}{=}(-1)^{|\phi|(|p|+|\Delta
|)+|\omega||p|}(\phi\circ\Delta)(\omega p),\quad\omega\in\mathcal{C}_{\star
}\Lambda_{k-1}(\mathcal{E}).
\]
Tensor products we need now on are over the algebra $\mathcal{C}_{\star}%
\Lambda_{k-1}(\mathcal{E})$ and  we shall simplify the notation by using
$\otimes$ instead of $\otimes_{\mathcal{C}_{\star}%
\Lambda_{k-1}(\mathcal{E})}.$

\begin{proposition}
There exists a unique $\mathcal{C}_{\star}\Lambda_{k-1}(\mathcal{E})$--module
homomorphism
\[
\mathcal{C}\mathrm{Diff}(P_{1},P_{2})\ni\Delta\longrightarrow\Delta^{Q}%
\in\mathcal{C}\mathrm{Diff}(P_{1}\otimes Q,P_{2}\otimes Q)
\]
such that
\[
\Delta^{Q}(p\otimes q,\phi\otimes\psi)(\omega)=(-1)^{|\psi|(|q|+|p|+|\Delta
|+|\phi|)+|q|(|\phi|+|\omega|)}\psi(\Delta(p,\phi)^{Q}(\omega q)),
\]
$p\in P_{1}$, $q\in Q$, $\phi\in P_{2}^{\ast}$, $\psi\in Q^{\ast}$, $\omega
\in\mathcal{C}_{\star}\Lambda_{k-1}(\mathcal{E})$.
\end{proposition}

Now, let $P_{3}$ be another
$\mathcal{C}_{\star}\Lambda_{k-1}(\mathcal{E})$--module. Then
\[
(\Delta_{2}\circ\Delta_{1})^{Q}=\Delta_{2}^{Q}\circ\Delta_{1}^{Q}
\]
for any $\Delta_{1}\in\mathcal{C}\mathrm{Diff}(P_{1},P_{2}), \;\Delta_{2}%
\in\mathcal{C}\mathrm{Diff}(P_{2},P_{3})$. This simple fact allows
to associate with a complex of $\mathcal{C}_{\star}%
\Lambda_{k-1}(\mathcal{E})$--modules
\[
\cdots\longrightarrow P_{i-1}\overset{\Delta_{i-1}}{\longrightarrow}%
P_{i}\overset{\Delta_{i}}{\longrightarrow}P_{i+1}\overset{\Delta_{i+1}%
}{\longrightarrow}\cdots
\]
connected by $\mathcal{C}$--differential operators $\Delta_i$'s a
new one, namely,
\[
\cdots\longrightarrow P_{i-1}\otimes Q\overset{\Delta_{i-1}^{Q}}%
{\longrightarrow}P_{i}\otimes Q\overset{\Delta_{i}^{Q}}{\longrightarrow
}P_{i+1}\otimes Q\overset{\Delta_{i+1}^{Q}}{\longrightarrow}\cdots
\]

\section{$(k-1)$--IDF--symmetries and $\Lambda_{k-1}\mathcal{C}$--Spectral
Sequence of a System of PDEs}

Let $\mathcal{E}$ be as above. The universal linearization
$\ell_{\Lambda_{k-1}F}^{\{k\}}:\Lambda_{k-1}\varkappa\longrightarrow
\Lambda_{k-1}P$ of $\Lambda_{k-1}F$ is a
$\mathcal{C}$--differential operator and as such (proposition
\ref{PropRestr}) restricts to $\mathcal{E}$.
\begin{theorem}
\label{ThSym}There is a Lie--algebra isomorphism $\Lambda_{k-1}\mathrm{Sym}%
(\mathcal{E})\simeq\ker(\ell_{\Lambda_{k-1}F}^{{\{k\}}}|_{\mathcal{E}}%
)\subset\Lambda_{k-1}\varkappa |_{\mathcal{E}}$.
\end{theorem}

An exact description of this isomorphism is as follows. Let $\overline{\chi}%
\in\ker(\ell_{\Lambda_{k-1}F}^{\{k\}}|_{\mathcal{E}})\subset\Lambda
_{k-1}\varkappa|_{\mathcal{E}}$.
Then $\overline{\chi}=R_{\mathcal{E}%
}^{\Lambda_{k-1}\varkappa}(\chi)$ for some
$\chi\in\Lambda_{k-1}\varkappa $. The evolutionary derivation
$\rE_{\chi}:\Lambda_{k-1}\longrightarrow\Lambda_{k-1}$ can be
restricted to $\mathcal{E}$, i.e., there exists a unique derivation
$\rE_{\chi}|_{\mathcal{E}}%
\in\mathrm{D}(\Lambda_{k-1}(\mathcal{E}),\Lambda_{k-1}(\mathcal{E}))$
such that $i_{\mathcal{E}}^{\ast}\circ\rE_{\chi}=\rE_{\chi
}|_{\mathcal{E}}\circ i_{\mathcal{E}}^{\ast}$. Moreover,
$\rE_{\chi }|_{\mathcal{E}}\in $
$\mathrm{D}_{\mathcal{C}}(\Lambda_{k-1}(\mathcal{E}))$. Then the
derivation that corresponds to $\overline{\chi}$ via the
isomorphism of theorem
\ref{ThSym} is $[\rE_{\chi}|_{\mathcal{E}%
}]_{\mathcal{C}\mathrm{D}(\Lambda_{k-1}(\mathcal{E}))}\in\Lambda
_{k-1}\mathrm{Sym}(\mathcal{E})$.

Put
\[
\Lambda_{k-1}Q\overset{\mathrm{def}}{=}\{\nabla|_{\mathcal{E}}\circ
\ell_{\Lambda_{k-1}F}^{\{k\}}|_{\mathcal{E}}\;|\;\nabla\in\mathcal{C}\mathrm{Diff}(\Lambda
_{k-1}P,\Lambda_{k-1})\}\subset\mathcal{C}\mathrm{Diff}(\Lambda_{k-1}%
\varkappa|_{\mathcal{E}},\Lambda_{k-1}(\mathcal{E}))
\]
and note that
$\mathcal{C}\mathrm{Diff}_{(p-1)}(\Lambda_{k-1}\varkappa
|_{\mathcal{E}};\Lambda_{k-1}Q)$ is naturally embedded into $\mathcal{C}%
\mathrm{Diff}_{(p)}(\Lambda_{k-1}\varkappa|_{\mathcal{E}};\Lambda
_{k-1}(\mathcal{E}))$. In view of this embedding the following
definition
\[
\Lambda_{k-1}Q_{p}\overset{\mathrm{def}}{=}\mathrm{alt}_{p}(\mathcal{C}%
\mathrm{Diff}_{(p-1)}(\Lambda_{k-1}\varkappa|_{\mathcal{E}};\Lambda
_{k-1}Q))\subset\mathcal{C}\mathrm{Diff}_{(p)}^{\mathrm{alt}}(\Lambda
_{k-1}\varkappa|_{\mathcal{E}};\Lambda_{k-1}(\mathcal{E}))
\]
makes sense. In particular, $\Lambda_{k-1}Q_{1}=\Lambda_{k-1}Q$.

\begin{proposition}
There is a $\Lambda_{k-1}(\mathcal{E})$--module isomorphism
\[
\mathcal{C}^{p}\Lambda_{k}^{p}(\mathcal{E})\simeq\mathcal{C}\mathrm{Diff}%
_{(p)}^{\mathrm{alt}}(\Lambda_{k-1}\varkappa|_{\mathcal{E}};\Lambda
_{k-1}(\mathcal{E}))/\Lambda_{k-1}Q_{p}.
\]

\end{proposition}

According to this proposition the homomorphism
\[
\lambda_{\Lambda_{k-1}F}:\mathcal{C}\mathrm{Diff}(\Lambda_{k-1}P|_{\mathcal{E}%
};\Lambda_{k-1}(\mathcal{E}))\longrightarrow\mathcal{C}\mathrm{Diff}%
(\Lambda_{k-1}\varkappa|_{\mathcal{E}};\Lambda_{k-1}(\mathcal{E})),
\]
sending  $\nabla\in\mathcal{C}%
\mathrm{Diff}(\Lambda_{k-1}P|_{\mathcal{E}};\Lambda_{k-1}(\mathcal{E}))$
to $\nabla\circ\ell_{\Lambda_{k-1}F}^{\{k\}}$, makes the following
sequence of $\Lambda
_{k-1}(\mathcal{E})$--homomorphisms%
\[
\mathcal{C}\mathrm{Diff}(\Lambda_{k-1}P|_{\mathcal{E}};\Lambda_{k-1}%
(\mathcal{E}))\overset{\lambda_{\Lambda_{k-1}F}}{\longrightarrow}\mathcal{C}%
\mathrm{Diff}(\Lambda_{k-1}\varkappa|_{\mathcal{E}};\Lambda_{k-1}(\mathcal{E}%
))\longrightarrow\mathcal{C}\Lambda_{k}^{1}(\mathcal{E})\longrightarrow
0,
\]
exact.

Recall (see \cite{b99,v01}) that the equation $\mathcal{E}_{F}$ is
said $\ell$--\emph{normal} if the sequence of
$C^{\infty}(\mathcal{E})$--homomorphisms
\[
0\longrightarrow\mathcal{C}\mathrm{Diff}(P|_{\mathcal{E}};C^{\infty
}(\mathcal{E}))\overset{\lambda_{F}}{\longrightarrow}\mathcal{C}%
\mathrm{Diff}(\varkappa|_{\mathcal{E}};C^{\infty}(\mathcal{E}))\longrightarrow
\mathcal{C}\Lambda^{1}(\mathcal{E})\longrightarrow0
\]
with $\lambda_{F}(\nabla)\overset{\mathrm{def}}{=}\nabla\circ
\ell_{F}$,
$\nabla\in\mathcal{C}\mathrm{Diff}(P|_{\mathcal{E}};C^{\infty
}(\mathcal{E}))$, is exact. The \textquotedblleft iterated\textquotedblright analogue of this notion
is as follows.

\begin{definition}
An equation $\mathcal{E}_{F}$ is said
$\ell^{\{k\}}$--\emph{normal} if the sequence of
$\Lambda_{k-1}(\mathcal{E})$--homomorphisms
\[
0\longrightarrow\mathcal{C}\mathrm{Diff}(\Lambda_{k-1}P|_{\mathcal{E}}%
;\Lambda_{k-1}(\mathcal{E}))\overset{\lambda_{\Lambda_{k-1}F}}{\longrightarrow
}\mathcal{C}\mathrm{Diff}(\Lambda_{k-1}\varkappa|_{\mathcal{E}};\Lambda
_{k-1}(\mathcal{E}))\longrightarrow\mathcal{C}\Lambda_{k}^{1}(\mathcal{E}%
)\longrightarrow0
\]
is exact, i.e., $\ker\lambda_{\Lambda_{k}F}=0$.
\end{definition}

\begin{proposition}
If $\mathcal{E}$ is $\ell$--normal, then it is
$\ell^{(\{k\}}$--normal as well for any $k>1$.
\end{proposition}

The zeroth column differential
$d_{0}^{0,\bullet}:\mathcal{H}\Lambda_{k}(\mathcal{E})\longrightarrow
\mathcal{H}\Lambda_{k}(\mathcal{E})$   in the zeroth term of the
$\Lambda_{k-1}\mathcal{C}$--spectral sequence of $\mathcal{E}$ is
a $\mathcal{C}$--differential operator. Consider the operator
$[d_{0}^{0,\bullet}]_{p}\overset{\mathrm{def}}{=}(d_{0}^{0,\bullet
})^{\mathcal{C}_{\star}\Lambda_{k}^{1}(\mathcal{E})\otimes\mathcal{C}_{\star
}\Lambda_{k}^{p-1}(\mathcal{E})}$,
\[
\lbrack d_{0}^{0,\bullet}]_{p}:\mathcal{C}_{\star}\Lambda_{k}^{1}%
(\mathcal{E})\otimes\mathcal{C}_{\star}\Lambda_{k}^{p-1}(\mathcal{E}%
)\otimes\mathcal{H}\Lambda_{k}(\mathcal{E})\longrightarrow\mathcal{C}_{\star
}\Lambda_{k}^{1}(\mathcal{E})\otimes\mathcal{C}_{\star}\Lambda_{k}%
^{p-1}(\mathcal{E})\otimes\mathcal{H}\Lambda_{k}(\mathcal{E})
\]
and denote by
\begin{equation}
\mathrm{alt}_{p}:\mathcal{C}_{\star}\Lambda_{k}^{1}(\mathcal{E})\otimes
\mathcal{C}_{\star}\Lambda_{k}^{p-1}(\mathcal{E})\otimes\mathcal{H}\Lambda
_{k}(\mathcal{E})\longrightarrow\mathcal{C}_{\star}\Lambda_{k}^{p}%
(\mathcal{E})\otimes\mathcal{H}\Lambda_{k}(\mathcal{E}) \label{proj}%
\end{equation}
the alternation map. Then a natural isomorphism
$\mathcal{C}_{\star
}\Lambda_{k}^{p}(\mathcal{E})\otimes\mathcal{H}\Lambda_{k}(\mathcal{E}%
)\simeq\Lambda_{k-1}\mathcal{C}E_{0}^{p,\bullet}$ takes place.
Moreover, the inclusion
$\mathcal{C}_{\star}\Lambda_{k}^{p}(\mathcal{E})\otimes\mathcal{H}\Lambda
_{k}(\mathcal{E})\subset\mathcal{C}_{\star}\Lambda_{k}^{1}(\mathcal{E}%
)\otimes\mathcal{C}_{\star}\Lambda_{k}^{p-1}(\mathcal{E})\otimes
\mathcal{H}\Lambda_{k}(\mathcal{E})$ is a right inverse of $\mathrm{alt}_{p}$
so that $\Lambda_{k-1}\mathcal{C}E_{0}^{p,\bullet}$ is a direct summand in
$\mathcal{C}_{\star}\Lambda_{k}^{1}(\mathcal{E})\otimes\mathcal{C}_{\star
}\Lambda_{k}^{p-1}(\mathcal{E})\otimes\mathcal{H}\Lambda_{k}(\mathcal{E})$.
More precisely, the complex $(\Lambda_{k-1}\mathcal{C}E_{0}^{p,\bullet}%
,d_{0}^{p,\bullet})$ is a direct summand in the complex $(\mathcal{C}_{\star
}\Lambda_{k}^{1}(\mathcal{E})\otimes\mathcal{C}_{\star}\Lambda_{k}%
^{p-1}(\mathcal{E})\otimes\mathcal{H}\Lambda_{k}(\mathcal{E}),[d_{0}%
^{0,\bullet}]_{p})$.

Now consider the restriction to $\mathcal{E}$ of the adjoint
operator of $\ell_{\Lambda_{k-1}F}^{\{k\}}$:
\[
\widehat{\ell}_{\Lambda_{k-1}F}^{\{k\}}|_{\mathcal{E}}:\widehat{\Lambda
_{k-1}P}|_{\mathcal{E}}\longrightarrow\widehat{\Lambda_{k-1}\varkappa
}|_{\mathcal{E}}%
\]
and its multiple extension
\[
{}[\widehat{\ell}_{\Lambda_{k-1}F}^{\{k\}}|_{\mathcal{E}}]_{p}:\widehat
{\Lambda_{k-1}P}|_{\mathcal{E}}\otimes\mathcal{C}_{\star}\Lambda_{k}%
^{p}(\mathcal{E})\longrightarrow\widehat{\Lambda_{k-1}\varkappa}%
|_{\mathcal{E}}\otimes\mathcal{C}_{\star}\Lambda_{k}^{p}(\mathcal{E})
\]
defined as ${}[\widehat{\ell}_{\Lambda_{k-1}F}^{\{k\}}|_{\mathcal{E}}%
]_{p}\overset{\mathrm{def}}{=}(\widehat{\ell}_{\Lambda_{k-1}F}^{\{k\}}%
|_{\mathcal{E}})^{\mathcal{C}_{\star}\Lambda_{k}^{p}(\mathcal{E})}$.

\begin{proposition}
\quad

\begin{itemize}
\item $H^{q}(\mathcal{C}_{\star}\Lambda_{k}^{1}(\mathcal{E})\otimes
\mathcal{C}_{\star}\Lambda_{k}^{p-1}(\mathcal{E})\otimes\mathcal{H}\Lambda
_{k}(\mathcal{E}),[d_{0}^{0,\bullet}]_{p})\simeq0$ for $q<n-1$.

\item $H^{n-1}(\mathcal{C}_{\star}\Lambda_{k}^{1}(\mathcal{E})\otimes
\mathcal{C}_{\star}\Lambda_{k}^{p-1}(\mathcal{E})\otimes\mathcal{H}\Lambda
_{k}(\mathcal{E}),[d_{0}^{0,\bullet}]_{p})\simeq\ker{}[\widehat{\ell}%
_{\Lambda_{k-1}F}^{\{k\}}|_{\mathcal{E}}]_{p-1}$.

\item $H^{n}(\mathcal{C}_{\star}\Lambda_{k}^{1}(\mathcal{E})\otimes
\mathcal{C}_{\star}\Lambda_{k}^{p-1}(\mathcal{E})\otimes\mathcal{H}\Lambda
_{k}(\mathcal{E}),[d_{0}^{0,\bullet}]_{p})\simeq\operatorname{coker}%
{}[\widehat{\ell}_{\Lambda_{k-1}F}^{\{k\}}|_{\mathcal{E}}]_{p-1}$.
\end{itemize}
\end{proposition}

In the following theorem that collects the above results the map
$\mathrm{alt}_{p}$, abusing the notation, stands for the induced
by projection (\ref{proj}) map in cohomology.
\begin{theorem}
[Two Lines Theorem]Let $\mathcal{E}=\mathcal{E}_{F}^{\infty}$ be
the infinite prolongation of an $\ell$--normal equation
$\mathcal{E}_{F}=\{F=0\}$. Then $\mathcal{E}$ is
$\ell^{\{k\}}$--normal and

\begin{itemize}
\item $\Lambda_{k-1}\mathcal{C}E_{1}^{0,\bullet}(\mathcal{E})\simeq
\Lambda_{k-2}\mathcal{C}E_{1}^{\bullet,\bullet}(\mathcal{E})$.

\item $\Lambda_{k-1}\mathcal{C}E_{1}^{p,q}(\mathcal{E})=0$ for any $p>0$ and
$q<n-1$.

\item $\Lambda_{k-1}\mathcal{C}E_{1}^{p,n-1}(\mathcal{E})\simeq\mathrm{alt}%
_{p}(\ker{}[\widehat{\ell}_{\Lambda_{k-1}F}^{\{k\}}|_{\mathcal{E}}]_{p-1})$
for any $p>0$.

\item $\Lambda_{k-1}\mathcal{C}E_{1}^{p,n}(\mathcal{E})\simeq\mathrm{alt}%
_{p}(\operatorname{coker}{}[\widehat{\ell}_{\Lambda_{k-1}F}^{\{k\}}%
|_{\mathcal{E}}]_{p-1})$ for any $p>0$.
\end{itemize}
\end{theorem}

\section{Short Conclusion}

Limits of a short note do not allow us neither to be completely
conceptual in the exposition, nor to present all obtained results
in full generality. Both require introducing numerous new notions
and constructions that can be done only in frames of the
forthcoming detailed and systematical exposition. The same
concerns various applications, even quite direct ones as, for
instance, secondary versions of notes \cite{vv06,vv06b,vv07c}. By concluding we state that the \textquotedblleft iterated\textquotedblright analogues
of the $l$--lines theorem and all results contained in the last
chapter of book \cite{v01} take place.


\begin{thebibliography}{9}                                                                                                %


\bibitem {vv07}A.~M.~Vinogradov, L.~Vitagliano, \emph{Dokl.~Math.}, \emph{to
appear in}, see also The Diffiety Inst.~Preprint Series, DIPS 5/06 or arXive: math.DG/0610917.

\bibitem {b99}A.~V.~Bocharov, V.~N.~Chetverikov, S.~V.~Duzhin,
N.~G.~Khor'kova, I.~S.~Krasil'shchik, A.~V.~Samokhin, Yu.~N.~Torkhov,
A.~M.~Verbovetsky, and A.~M.~Vinogradov, \emph{Symmetries and Conservation
Laws for Differential Equations of Mathematical Physics}, AMS
\textquotedblleft Translation of Mathematical Monographs\textquotedblright%
\ vol. \textbf{182} (1999).

\bibitem {vv06}A.~M.~Vinogradov, L.~Vitagliano, \emph{Dokl.~Math.~}%
\textbf{73}, n$^{\circ}$ 2 (2006) 169, see also The Diffiety Inst.~Preprint
Series, DIPS 1/06 or arXive: math.DG/0605113.

\bibitem {vv07b}A.~M.~Vinogradov, L.~Vitagliano, \emph{Dokl.~Math.}, \emph{to
appear in}, see also The Diffiety Inst.~Preprint Series, DIPS 1/07 or arXive: math.DG/0703661.

\bibitem {v01}A.~M.~Vinogradov, \emph{Cohomological Analysis of Partial
Differential Equations and Secondary Calculus}, AMS \textquotedblleft
Translation of Mathematical Monographs\textquotedblright\ vol. \textbf{204} (2001).

\bibitem {kv98}I.~S.~Krasil'shchik, A.~M.~Verbovetsky, \emph{Homological
Methods in Equation of Mathematical Physics} Open Education (Opava) 1998. See
also The Diffiety Inst.~Preprint Series, DIPS 7/98,
\texttt{http://diffiety.ac.ru/preprint/98/08\_98abs.htm}.

\bibitem {vv06b}A.~M.~Vinogradov, L.~Vitagliano, \emph{Dokl.~Math.}~\textbf{73}, n$^{\circ}$ 2 (2006) 182, see also The Diffiety Inst.~Preprint
Series, DIPS 3/06 or arXive: math.DG/06009287.

\bibitem {vv07c}A.~M.~Vinogradov, L.~Vitagliano, \emph{Dokl.~Math.}~\textbf{75}, \emph{to appear in}, see also The Diffiety Inst.~Preprint
Series, DIPS 4/06 or arXiv: math.DG/0610914..

\end{thebibliography}
\end{document}